\documentclass[12pt]{article}
\usepackage{amsmath, amssymb, amsthm, enumitem} 

\usepackage[T1]{fontenc}

\newtheorem{theorem}{Theorem}

\newtheorem{proposition}{Proposition}

\newtheorem{definition}{Definition}

\newtheorem{conjecture}{Conjecture}

\DeclareMathOperator{\diam}{diam}
\DeclareMathOperator{\id}{id}
\DeclareMathOperator{\Int}{Int}
\DeclareMathOperator{\Ext}{Ext}
\DeclareMathOperator{\Emb}{Emb}

\begin{document}

\title{An answer to a question of J.W.~Cannon and S.G.~Wayment
}
\author{Olga Frolkina\\
Chair of General Topology and Geometry\\
Faculty of Mechanics and Mathematics\\
M.V.~Lomonosov Moscow State University\\
and Moscow Center for Fundamental and Applied Mathematics\\
Leninskie Gory 1, GSP-1,
Moscow 119991, Russia\\
E-mail: olga-frolkina@yandex.ru
}

\date{}

\maketitle

\renewcommand{\thefootnote}{}



\renewcommand{\thefootnote}{\arabic{footnote}}
\setcounter{footnote}{0}

\begin{abstract}
Solving R.J.~Daverman's problem,
V.~Krushkal described 
sticky Cantor sets in
$\mathbb R^N$ for $N\geqslant 4$;
these sets cannot be isotoped off of itself by small ambient isotopies.
Using Krushkal sets, we answer a question
of J.W.~Cannon and S.G.~Wayment (1970).
Namely, for $N\geqslant 4$ we construct compacta $X\subset \mathbb R^N$
with the following two properties:
some sequence $\{ X_i \subset \mathbb R^N \setminus X, \ i\in\mathbb N \}$ 
converges homeomorphically to $X$, but 
there is no uncountable family of pairwise disjoint
sets $Y_\alpha \subset \mathbb R^N$
each of which is embedded equivalently to~$X$.
\end{abstract}

\section{Introduction}

For any non-empty compactum $\mathfrak X$,
the space $\Emb (\mathfrak X,\mathbb R^N)$
of embeddings
$\mathfrak X\hookrightarrow \mathbb R^N$
endowed with the metric
$\rho (f,g) = \sup \{ d(fx,gx) \ | \ x\in\mathfrak X \}$
is separable.
An immediate consequence is the following known result:
Let $\mathcal F$ be an uncountable family of
mutually disjoint homeomorphic copies of
a compactum $\mathfrak X$
in $\mathbb R^N$;
then
there exists a sequence $X_0,X_1,X_2,\ldots \in \mathcal F$ 
such that 
$X_i \neq X_0$ (hence
$X_i \cap X_0 = \emptyset $) for each $i>0$, and
$\{ X_i, \ i\in\mathbb N \}$ converges homeomorphically to~$X_0$.

\begin{definition}\label{convergence}
A sequence 
of sets  $X_1,X_2,\ldots \subset \mathbb R^N$
converges homeo\-morphi\-cal\-ly to a set
$X_0\subset \mathbb R^N $
if
for each
$\varepsilon > 0$ there is an integer~$n$ such that $i\geqslant n$ implies
the existence of a homeomorphism
$h_i: X_0\cong X_i$
which moves no point more than $\varepsilon $.
\end{definition}

In 1970, J.W.~Cannon and S.G.~Wayment raised the question \cite{Cannon-Wayment}:

Suppose that $X_0, X_1,X_2,\ldots \subset \mathbb R^N$
is a sequence of pairwise disjoint
continua that converges homeomorphically to $X_0$.
Is it true that 
there exists an uncountable family 
of pairwise disjoint
homeomorphic copies of $X_0$ in $\mathbb R^N$~?
(According to the result obtained in 1993 by
E.K.~van Douwen,
an equivalent question is whether
$X_0\times\mathcal C$ embeds in $\mathbb R^N$,
where $\mathcal C$ is the Cantor set \cite{vanDouwen}. 
H.~Becker, F.~van Engelen and J.~van Mill
give a short proof of van Douwen's theorem, together 
with a generalization and an example
of a $\sigma $-compact space which contains an uncountable
family of pairwise disjoint copies of the circle $S^1$
but not $S^1 \times \mathcal C$ \cite{BEM}.
S.~Todor\v{c}evi\'c
gives a combinatorial characterization for
the embeddability 
of $X\times \mathcal C$ into a complete separable space $Y$;
as a corollary, $X\times \mathcal C$ embeds into $Y$ iff
$X\times \mathbb Q$ does \cite{Todorcevic}.)

This question remains open. 
Cannon and Wayment
obtained a positive answer under the additional 
assumption: $X_0,X_1,X_2,\ldots $
are embedded in $\mathbb R^N$ equivalently to each other
\cite[Thm. 1]{Cannon-Wayment}.

\begin{definition} 
Two subsets $X,X'\subset \mathbb R^N$ 
are \emph{ambiently homeomorphic}
(or \emph{equivalently embedded})
if there exists a homeomorphism $h$
of $\mathbb R^N$ onto itself such that
$h(X)=X'$.
In this case, we write
$h: (\mathbb R^N,X)\cong (\mathbb R^N,X')$.
\end{definition}

But 
even under this assumption, 
it is not always possible to find a desired uncountable family
so that all of its elements were
embedded 
equivalently to the limit
set $X_0\subset \mathbb R^N$.
This is confirmed by examples.

For $N=2$, see \cite[Examples 3, 4]{Bing-snakelike},
\cite[Example 2]{Roberts};
the spaces $X_0$ have rather complicated structure.

For $N=3$ or $N\geqslant 5$,
Cannon and Wayment constructed
a sequence 
$X_0,X_1,X_2,\ldots $ 
of pairwise disjoint wild $(N-1)$-spheres
in $\mathbb R^N$ such that:
$\{ X_i ,\ i\in\mathbb N\}$ converges homeomorphically to $X_0$;
$(\mathbb R^N, X_i)\cong (\mathbb R^N,X_0)$
for each $i$;
but there is no possibility
to find uncountably many pairwise
disjoint $(N-1)$-spheres in $\mathbb R^N$
each of which is embedded equivalently
to $X_0$
\cite[p.~568--570]{Cannon-Wayment}.

Here it is appropriate to recall

\begin{definition} 
A subset $X\subset \mathbb R^N$ 
homeomorphic to the $(N-1)$-sphere
is called \emph{flat}
if it is ambiently homeomorphic
to the standard unit sphere,
and \emph{wild} otherwise.
\end{definition}

The case of 
$N=4$ was not covered by the examples of
J.W.~Cannon and S.G.~Wayment, and remained an open problem 
\cite[p. 569--570]{Cannon-Wayment}.
In the present paper, 
assuming $N\geqslant 4$,
we construct new series of embeddings:
for $(N-1)$-spheres (Theorem \ref{CW-answer-sphere});
for a wider class of compacta of positive dimension 
(Theorem \ref{CW-answer-gen});
and for Cantor sets (Theorem \ref{Cantor}).
(Cannon and Wayment restricted their question
by connected sets, but the case of Cantor sets is also non-trivial.)

To some extent, the possibility to find
the desired sequences
is explained by the fact that
a small homeomorphism $X_0\cong X_i$ may have
no small extending homeomorphism 
$\mathbb R^N\cong\mathbb R^N$.
Compare with \cite[Thm.~4, Cor.~3,~4]{Frolkina2022}.

For a review of notions and results concerning 
wild embeddings, refer to
\cite{Keldysh},
\cite{Rushing}, \cite{DV},
\cite{BC}, \cite{Burgess-75}, \cite{Daverman}, \cite{Cernavsky}.
A thorough discussion of famous 
three-dimensional examples belonging to L.~Antoine, J.~Alexander,
R.H.~Fox and E.~Artin can also be found in
\cite{Moise}.

\subsection{Examples of Cannon and Wayment: the idea}

For $N=3$, each of $X_0,X_1,X_2,\ldots $  is embedded 
equivalently to
the Fox--Artin sphere
\cite[Thm.~20.5]{Moise},
\cite[2.8.3]{DV} which is wild at only one point.
It is impossible to find uncountably many
mutually disjoint Fox-Artin spheres in $\mathbb R^3$
by the result of
R.H.~Bing: 
one can not place in
$\mathbb R^3$
uncountably many
pairwise disjoint wild closed
surfaces.
(Scheme of proof is given in
\cite{Bing-abstr};
full proof applies results of
\cite{Bing59} and \cite{Bing-TAMS61}.
Idea of Bing's proof can also be found in
\cite[Thm.~3.6.1]{BC}.)

Suppose that $N\geqslant 5$.
Again, the sequence 
$X_0,X_1,X_2,\ldots $ consists of wild
$(N-1)$-spheres.
But there is no analogue
of Fox-Artin sphere:
if $N\geqslant 4$
and an $(N-1)$-sphere $\Sigma \subset \mathbb R^N$ is not flat,
then the set of points at which $\Sigma $ fails to be locally flat
contains more than one point
\cite{Cantrell},
\cite[Thm. 3.3.1]{Rushing}, 
\cite[Thm. 2.9.3]{DV}
(a stronger ``local'' result is proved in
\cite{Cernavsky1966},
\cite{Cernavsky1968a},
\cite{Cernavsky2006},
\cite{Kirby2}; see also \cite[3A]{Daverman} for a review).
A desired example can be obtained analogously to the
Antoine's wild $2$-sphere
\cite[Examples 2.1.7, 4.7.1]{DV},
\cite[Thm. 18.7]{Moise},
the wildness points form a Cantor set
which is wild as a subset of
$\mathbb R^N$.
For $N\geqslant 5$
it is impossible
to find uncountably many pairwise
disjoint wild $(N-1)$-spheres in $\mathbb R^N$
\cite[Thm. 1, 2]{Bryant}, 
\cite[Thm. 10.5]{Burgess-75}, \cite[p.~383, Thm. 3C.2]{Daverman}.

The proof of the Bing impossibility theorem 
includes careful geometric analysis.
Its analogue for dimensions $\geqslant 5$ is based on such
deep results as
M.~Brown's version of the 
Generalized Sch{\"o}nflies Theorem
\cite{MBrown},
\cite[Cor. 5.1.3]{Keldysh},
\cite[Cor. 2.4.13]{DV} and 
the homotopy-theoretic criterion for local flatness
\cite{Cernavsky-criterion},
\cite{Daverman-criterion},
\cite[3B, 3C]{Daverman}, \cite[Thm. 7.6.1, 7.6.5]{DV}.

\subsection{How do our examples differ?}

For compacta of positive dimension,
(Theorems \ref{CW-answer-sphere},
\ref{CW-answer-gen})
our idea resembles that
of Cannon--Wayment. 
The key difference is the use
of sticky Cantor sets 
described by V.~Krushkal.
Each of our compacta contains a Krushkal Cantor set;
it is impossible to find uncountably many 
pairwise disjoint Krushkal Cantor sets in $\mathbb R^N$ for
$N\geqslant 4$;
hence it is impossible to find uncountably many 
pairwise disjoint homeomorphic copies of our compacta. 
This approach allows us to
cover all dimensions $N\geqslant 4$ and to
go beyond the case
of $(N-1)$-spheres.
This also 
eliminates the need in difficult 
Bing-type 
impossibility 
results.

\subsection{Notation and agreements}

$\mathbb R^N$ is the $N$-dimensional Euclidean space
with the standard metric~$d$. 
$O$ is its zero point.
$S^{N-1}$ is the standard unit sphere in $\mathbb R^{N}$.
For $a\in\mathbb R^N$ and $r>0$
let
$
U(a,r) = \{ x\in\mathbb R^N \ | \ d(x,a) < r \}
$
and 
$
B(a,r) = \{ x\in\mathbb R^N \ | \ d(x,a) \leqslant r \}
$
.
An $N$-ball in $\mathbb R^N$
is a set of the form $B(a,r)$.
But an $(N-1)$-sphere in $\mathbb R^N$
is any subset homeomorphic to $S^{N-1}$.

For any (possibly wild) $(N-1)$-sphere $\Sigma \subset \mathbb R^N$,
the complement
$\mathbb R^N\setminus \Sigma $
has two connected components
 (Jordan--Brouwer Separation Theorem).
The bounded component is denoted by
$\Int \Sigma $,
and the unbounded one --- by 
$\Ext \Sigma $.

$\partial  A$ is the boundary of a set $A\subset \mathbb R^N$.

The symbol $X \cong X'$ stands for a homeomorphism between $X$ and $X'$.
We write
$(\mathbb R^N , X)\cong (\mathbb R^N , X')$ 
if  there exists a homeomorphism of pairs, i.e.
$h: \mathbb R^N \cong \mathbb R^N$ such that
$h(X)= X'$.

For two subsets $X,X'\subset \mathbb R^N$,
a homeomorphism $h: X\cong X'$ is called
an $\varepsilon $-homeomorphism
if $d(x, h(x))\leqslant \varepsilon  $
for each $x\in X$.

As usual, $\id $ is the identity map.

Disjunctive family of sets is a family consisting of pairwise disjoint sets.

\section{Statements}

Let us first consider a particular case of $(N-1)$-spheres.

\begin{theorem}\label{CW-answer-sphere}
For each $N\geqslant 4$
there exists a sequence of pairwise disjoint
wild $(N-1)$-spheres
$\Sigma , X _1 , X_2, \ldots  \subset \mathbb R^N$
such that

1) 
$\{ X_i , \ i\in\mathbb N \}$
converges homeomorphically to $\Sigma $;

2) $(\mathbb R^N , X_i)  \cong (\mathbb R^N , \Sigma )$
for each $i\in\mathbb N$; and

3) any disjunctive family 
$\{ Y_\alpha \subset \mathbb R^N\}$
such that
$(\mathbb R^N , Y_\alpha ) \cong (\mathbb R^N , \Sigma )  $ for each $\alpha $,
contains no more than countably many elements.
\end{theorem}

More generally, we have

\begin{theorem}\label{CW-answer-gen}
Let $N\geqslant 4$.
Let $\mathfrak X \subset S^{N-1}$
be any compactum such that
for some point $q\in \mathfrak X$, some open neighbourhood
$U$ of $q$ in $S^{N-1}$, and some integer $k\in \{ 1,\ldots , N-1\}$
we have
$(U, U\cap \mathfrak X) \cong
(\mathbb R^{N-1} , \mathbb R^k)$.
Then there exist an embedding
$e: \mathfrak X \hookrightarrow \mathbb R^N$ 
and
a sequence
of pairwise disjoint compacta
$X _1 , X_2, \ldots  \subset \mathbb R^N \setminus e( \mathfrak X )$
with the properties:

1) $\{ X_i , i\in\mathbb N\}$ 
converges homeomorphically to $e( \mathfrak X )$;

2) $(\mathbb R^N , X_i)  \cong (\mathbb R^N , e( \mathfrak X ))$
for each $i\in\mathbb N$; and

3) 
if $\mathcal F$ is a disjunctive family 
of sets  $Y_\alpha \subset \mathbb R^N$
such that
$(\mathbb R^N , Y_\alpha ) \cong (\mathbb R^N , e( \mathfrak X ))  $ for each $\alpha $,
then $\mathcal F$ 
contains no more than countably many elements.
\end{theorem}

The case of a zero-dimensional set also has sense.

\begin{definition}
A subset $X\subset \mathbb R^N$ 
homeomorphic to the Cantor set
is called \emph{tame}
if it is ambiently homeomorphic
to the standard middle-thirds Cantor set
constructed in the straight line segment $[0,1]\times \{ 0\}^{N-1}\subset\mathbb R^N$.
Otherwise, $X$ is called \emph{wild}.
\end{definition}

Wild Cantor sets in $\mathbb R^N$ exist for each 
$N\geqslant 3$.
First wild Cantor sets
in $\mathbb R^3$ 
(now called Antoine's necklaces)
were described by L.~Antoine
in 1920--21
\cite[Sec. 18]{Moise}. He 
also proved that each Cantor set
in $\mathbb R^2$ is tame 
\cite[Thm. 13.7]{Moise}.
Another example of a wild Cantor set in $\mathbb R^3$
was constructed by P.S.~Urysohn in 1922--1923.

\begin{theorem}\label{Cantor}
For each $N\geqslant 4$ there exists a 
sequence of pairwise disjoint wild
Cantor sets $X, X _1 , X_2, \ldots \subset \mathbb R^N$
such that

1) 
$\{ X_i , i\in\mathbb N \}$ converges homeomorphically to $X$;

2)
$(\mathbb R^N , X_i)  \cong (\mathbb R^N , X )$
for each $i\in\mathbb N$; 
and

3) any disjunctive family 
of Cantor sets  $Y_\alpha \subset \mathbb R^N$
such that
$(\mathbb R^N , Y_\alpha ) \cong (\mathbb R^N , X)  $ for each $\alpha $,
is no more than countable.
\end{theorem}

\section{Proofs}

R.J.~Daverman conjectured
\cite[Conj. 1]{Daverman-question}
that
for any two Cantor sets $X$ and $X'$ in $\mathbb R^N$
and any
$\varepsilon >0$
there is an $\varepsilon $-homeomorphism $h:\mathbb R^N \cong \mathbb R^N$
such that
$X \cap h(X') = \emptyset $.
This is true for $N\leqslant 3$
\cite[Thm.~1]{Sher1969}.
However, for each $N\geqslant 4$
V.~Krushkal described
a sticky (wild) Cantor set in $\mathbb R^N$:
it cannot be
isotoped off of itself by any sufficiently small 
ambient isotopy \cite[Thm. 1.1]{Krushkal}. 
(Compare with \cite{Wright-pushing} and
\cite{Cannon-Wright}.)
Krushkal's construction 
suits our purpose because of

\begin{proposition}\label{K-sets}\cite[Cor. 3.5]{Frolkina2022}
Let $K\subset\mathbb R^N$ be a Krushkal Cantor set,
$N\geqslant 4$.
Suppose that
$\{ X_\alpha  \ | \ \alpha \in \mathcal A \}$ is a family 
of
pairwise disjoint Cantor sets in $\mathbb R^N$ such that
each $X_\alpha $ is
ambiently homeomorphic to $K$.
Then
$\mathcal A$ is no more than countable.
\end{proposition}

We will also use the following evident observation.

\begin{proposition}\label{equiv}
For a compactum $X_0\subset \mathbb R^N$,
the following are equivalent:

($\alpha $) there exists a sequence $X_1,X_2,\ldots $
of pairwise disjoint compacta in $\mathbb R^N\setminus X_0$ such that
$\{ X_i,\ i\in\mathbb N\}$ converges homeomorphically to $X_0$;

($\beta $) 
for any $\varepsilon >0$ there exists a homeomorphism
$h_{\varepsilon } : X_0 \hookrightarrow \mathbb R^N$
such that
$X_0 \cap h_{\varepsilon } (X_0) = \emptyset $
and 
$d(x, h_{\varepsilon } (x)) \leqslant \varepsilon $
for each $x\in X_0$.
\end{proposition}

\subsection{Proof of Theorem \ref{CW-answer-sphere}}

Denote
$$
W = [-1,1]^{N-1} \times [-2,0],
\quad
Q = \left[ -\frac12,\frac12\right] ^{N-1} \times [0,1],
\quad
P = W\cap Q .
$$

 Take an embedding $F:W\hookrightarrow \mathbb R^N$
 such that
 $F|_{\partial W \setminus P} = \id $
 and
 ${F(P) \subset Q}$.
 Let $\Sigma := F(\partial W)$.
We prove that any such $\Sigma $
satisfies 
condition ($\beta $)
of Proposition \ref{equiv}. 
This is done in four steps which
(as we can guess)
are implicit in \cite{Cannon-Wayment}.
We present a full proof
which 
will help the reader to verify the validity of Theorem \ref{CW-answer-gen}.
To get property 3), we will put additional restrictions on $\Sigma $;
this is done after Step 4.

{\bf Step 1.}
\emph{
For each $\varepsilon >0$ there exists an $\varepsilon $-homeomorphism
$\varphi _{\varepsilon } : \Sigma \hookrightarrow \mathbb R^N$
such that
$\varphi _{\varepsilon } (\Sigma ) \subset \Int \Sigma $
and 
$\varphi _{\varepsilon } (\Sigma )  $ is flat.
}

Denote $f:= F|_{\partial W}$, so that $\Sigma = f (\partial W)$.
For $\delta >0$, let $s_\delta : \mathbb R^N \to \mathbb R^N$
be the homothetic transformation which takes 
$(\xi _1 , \xi_2 , \ldots , \xi _N) \in \mathbb R^N$
to
$(\delta \xi _1 ,\delta \xi_2 , \ldots ,\delta \xi _N)$.
Put $\varphi _\varepsilon := f s_\delta f^{-1}$,
where $\delta $ is small enough, so that
$\varphi _\varepsilon $ 
moves each point a distance of at most $\varepsilon $.
Assuming $\delta < 1$, we get
$\varphi _\varepsilon  (\Sigma )\subset \Int \Sigma $.
By construction,
$\varphi _\varepsilon  (\Sigma ) $ is bicollared;
hence it is
flat \cite{MBrown}, 
\cite[Thm. 1.8.2]{Rushing}, \cite[Thm. 2.4.8]{DV}.

{\bf Step 2.}
\emph{
For each $\varepsilon >0$ there exists an $\varepsilon $-homeomorphism
$\psi _{\varepsilon } : \partial W \hookrightarrow \mathbb R^N$
such that
$\psi _{\varepsilon } (\partial W ) \subset \Ext (\partial W) $
and 
$(\mathbb R^N , \psi _{\varepsilon } (\partial W )  ) \cong (\mathbb R^N , \Sigma ) $.
}

For each $\lambda \in (0,1]$
define an embedding $f_\lambda : \partial W \hookrightarrow \mathbb R^N$
by
$$
f_\lambda (x) = 
\begin{cases}
 x &\text{for } x\in \partial W\setminus s_\lambda P, 
\\
s_\lambda f s_\lambda ^{-1} (x) &\text{for }
x\in s_\lambda P ,
\end{cases}
$$
where $s_\lambda $ is the homothetic transformation defined in Step~1.
In particular, $f_1 = f$.
For each $\lambda  \in (0,1]$,
we have $f_\lambda (\partial W)=(\partial W \setminus s_\lambda P) \cup s_\lambda f(P)$, and
it can be easily seen 
that 
$(\mathbb R^N,  f_\lambda (\partial W)) \cong
(\mathbb R^N , \Sigma )$.
 Finally,  take $\psi _\varepsilon := s_{1+\delta } f_\lambda $,
 where $\delta $ and $\lambda $ are sufficiently small
positive numbers.

{\bf Step 3.}
\emph{
For each 
flat $(N-1)$-sphere $\Gamma \subset \mathbb R^N$ and each
$\varepsilon >0$ there exists an $\varepsilon $-homeomorphism
$\kappa _{\varepsilon } :\Gamma \hookrightarrow \mathbb R^N$
such that
$\kappa _{\varepsilon } (\Gamma  ) \subset \Ext \Gamma  $
and 
$(\mathbb R^N , \kappa _{\varepsilon } (\Gamma  )  ) \cong (\mathbb R^N , \Sigma ) $.
}

Fix a homeomorphism 
$H: (\mathbb R^N , \Gamma ) \cong (\mathbb R^N , \partial W)$.
Put $\kappa _\varepsilon := H^{-1} \psi _\delta H$,
where $\psi _\delta $ is defined in Step~2, with
$\delta >0$ sufficiently small.

{\bf Step 4.}
\emph{
For each $\varepsilon >0$ there exists an $\varepsilon $-homeomorphism
$h _{\varepsilon } : \Sigma \hookrightarrow \mathbb R^N$
such that
$h _{\varepsilon } (\Sigma ) \subset \Int \Sigma $
and 
$(\mathbb R^N , h _{\varepsilon } (\Sigma )  ) \cong (\mathbb R^N , \Sigma ) $.
}

Take $\varphi _{\varepsilon /2} $ from Step~1.
Let $\gamma = \frac12 
\min \{ \varepsilon , d (\Sigma , \varphi _{\varepsilon /2} \Sigma )  \}$.
Applying Step~3 to $\varphi _{\varepsilon /2} \Sigma $ 
and $\gamma $,
construct a $\gamma $-homeomorphism
$\kappa _\gamma $ such that 
$\kappa _\gamma (\varphi _{\varepsilon /2} \Sigma )$ is embedded
equivalently to $\Sigma $.
The composition $h_\varepsilon := \kappa _\gamma \varphi _{\varepsilon /2} $
is the desired map.

Condition ($\beta $)
of Proposition \ref{equiv}
is thus verified.
Consequently
a corresponding
sequence $X_1,X_2,\ldots $
with properties 1) and 2) can be constructed.

To obtain property 3),
we put additional restrictions on $F$.
Namely, fix a Krushkal Cantor set
 $K\subset Q\setminus \partial Q$.
We require $F$ to be choosed so that $F(P) \supset K$. 
This
can be achieved
by the feelers technique which comes back to L.~Antoine,
see \cite[St. 4]{Frolkina2019}
for references.
Property 3) now follows from Proposition \ref{K-sets}.

\subsection{Proof of Theorem \ref{CW-answer-gen}}

Essentially the same as that for Theorem \ref{CW-answer-sphere}.
Identifying $S^{N-1}$ and $\partial W$,
we may assume that $\mathfrak X$ is embedded in $\partial W$
so that
$$\mathfrak X\cap P = 
 \left[ -\frac12,\frac12\right] ^{k} \times \{ 0 \}^{N-k} =: L_k 
 \quad
 \text{ and }
 \quad \mathfrak X\cap \widehat P = 
 \left[ -\frac34,\frac34\right] ^{k} \times \{ 0 \}^{N-k} , $$
 where
 $\widehat P =  \left[ -\frac34,\frac34\right] ^{N-1} \times \{ 0 \} $.
Fix a Krushkal Cantor set $K\subset Q\setminus \partial Q$.
Take an embedding
$F:W\hookrightarrow \mathbb R^N$
 such that
 $F|_{\partial W \setminus P} = \id $,
 $F(P) \subset Q $,
 and
 $K\subset F(L_1 )$.
 The desired embedding $e$ is defined as $F|_{\mathfrak X}$.

\subsection{Proof of Theorem \ref{Cantor}}

For each $m\in \mathbb N$
let
$$
a_m = \left(  \frac{1}{m}, 0 , 0, \ldots , 0 \right) \in\mathbb R^N
\quad
\text{ and }
\quad
r_m = \frac{1}{2(m+1)^2} .
$$
Closed balls $B_m := B(a_m,r_m)$, $m\in\mathbb N$,
are pairwise disjoint and
do not contain the zero point $O$.

Let $\tau_m: \mathbb R^N \cong \mathbb R^N$
be an orientation-preserving similarity transformation
that takes $B(O,1)$ onto $B\left( a_m , \frac{r_m}{2} \right)$.

Fix a Krushkal Cantor set $K\subset U(O,1) $.
Let $K_m := \tau _m (K)$ and
$X:= K_1 \cup K_2 \cup K_3 \cup \ldots \cup \{ O \}$.
By the Brouwer Characterization 
Theorem \cite[Thm. 12.8]{Moise}, $X$ is a Cantor set.

In the rest of the proof, we show that $X$ satisfies condition ($\beta $)
of Proposition \ref{equiv}.
Together with Proposition \ref{K-sets},
this proves that $X$ is the desired set.

Observe that
$\diam B_m \leqslant \varepsilon $
for $m>m_0 := \left[ \frac{1}{\sqrt{\varepsilon }} -1 \right] $.

For each $m>m_0$, take 
an $N$-ball
$D_m \subset U(a_m,r_m)
\setminus
B\left( a_m , \frac{r_m}{2} \right)$
and an orientation-preserving similarity transformation
$\sigma _m : \mathbb R^N \cong \mathbb R^N$
such that
$\sigma _m (B(O,1)) = D_m$.
Define $L_m := \sigma _m (K)$.
Fix any homeomorphism
$g_m : K_m \cong L_m$.
By construction, $g_m$ is an $\varepsilon $-homeomorphism.

Now let $m \in \{ 1 , 2, \ldots , m_0 \}$.
Decompose $K_m$
as the union
$K_m = K_{m,1} \cup K_{m,2}\cup\ldots \cup K_{m,i_m}$,
where
each $K_{m,j}$ is a Cantor set with $\diam K_{m,j} < \varepsilon $,
and
$K_{m,j} \cap K_{m,q} = \emptyset $ if $j\neq q$.
For each 
$j=1,2,\ldots , i_m$ take an $N$-ball
$D_{m,j}$ such that
$\diam (K_{m,j}\cup D_{m,j}) \leqslant \varepsilon $,
$D_{m,j} \cap X = \emptyset $,
and
$D_{m,j} \cap D_{m,q} = \emptyset $
if $j\neq q$.
Take an orientation-preserving similarity transformation
$\sigma _{m,j} : \mathbb R^N \cong \mathbb R^N$
such that
$\sigma _{m,j} (B(O,1)) = D_{m,j}$.
Denote $L_{m,j} := \sigma _{m,j} (K)$.
Fix a homeomorphism
$g_{m,j} : K_{m,j} \cong L_{m,j}$.
By construction, $g_{m,j}$ is an $\varepsilon $-homeomorphism.

Now put
\begin{multline*}
Y:= (L_{1,1} \cup L_{1,2} \cup \ldots \cup L_{1, i_1})\cup
\\
\cup  (L_{2,1} \cup L_{2,2} \cup \ldots \cup L_{2, i_2})
 \cup
 \ldots
\cup
 (L_{m_0 , 1} \cup \ldots \cup L_{m_0  , i_{m_0 }}) \cup 
 \\
 \cup L_{m_0+1} \cup L_{m_0+2}\cup L_{m_0+3}\cup
 \ldots
\cup \{ O \} .
\end{multline*}
By the Brouwer Characterization Theorem, $Y$ is a Cantor set.
By construction, $X\cap Y = \emptyset $.
Define a homeomorphism $h_\varepsilon : X\cong Y$ by
$$
h_\varepsilon  (x) = 
\begin{cases}
 O &\text{for } x=O, 
\\
g_m(x) &\text{for } x\in K_m,\ \  m>m_0,
\\
g_{m,j} (x) &\text{for } x\in K_{m,j}, \  1\leqslant m\leqslant m_0,
\ 1\leqslant j\leqslant i_m .
\end{cases}
$$
Clearly $h_\varepsilon $ is an $\varepsilon $-homeomorphism.

It remains to show that
$(\mathbb R^N ,X)\cong (\mathbb R^N , Y)$.
Order the balls that parti\-ci\-pate in the construction 
of $Y$ alphabetically
into the sequence
$$
Z_1:= D_{1,1} , \ \ 
Z_2:= D_{1,2} , \ \ 
\ldots , \ \ 
Z_{i_1}:= D_{1, i_1} ;
$$
$$
Z_{i_1 + 1 } := D_{2,1}, \ \ 
Z_{i_1 + 2 } := D_{2,2},\ \ 
\ldots  , 
 \ \ 
Z_{i_1 + i_2} := D_{2 , i_2 };
 $$
  $$
 Z_{i_1 + i_2 + \ldots + i_{m_0 - 1 } +1 } := D_{m_0 ,1}, \  \
\ldots  , \ \ 
Z_{i_1 + i_2 + \ldots + i_{m_0 - 1 }  + i_{m_0} } := D_{m_0  ,i_{m_0}};
$$
$$
\ldots 
$$
$$
Z_{i_1 + i_2 + \ldots + i_{m_0 - 1 }  + i_{m_0} + 1}:= D_{m_0 + 1},
$$
$$
Z_{i_1 + i_2 + \ldots + i_{m_0 - 1 }  + i_{m_0} + 2}:=D_{m_0 + 2},
$$
$$
Z_{i_1 + i_2 + \ldots + i_{m_0 - 1 }  + i_{m_0} + 3}:=D_{m_0 + 3},
$$
$$
\ldots
$$
Let $T$ be a parallel translation of $\mathbb R^N$ such that
$$
B_1\cup B_2 \cup B_3 \cup \ldots \cup \{ O \}
\quad
\text{ and }\
\quad
T(Z_1\cup Z_2 \cup Z_3\cup \ldots \cup \{ O \})
$$
can be separated by a round sphere
of the form $\partial (B(a,r))$.
In the rest of the proof, 
we construct a homeomorphism $G:\mathbb R^N \cong R^N$
such that $G(X) = T(Y)$.
(Consequently the homeomorphism 
$T^{-1}G : \mathbb R^N \cong \mathbb R^N$
takes $X$ onto $Y$, as desired.)

Construct a sequence of pairwise disjoint connected open sets
$V_\ell \subset \mathbb R^N$, $\ell \in\mathbb N$, such that
$B_\ell \cup T(Z_\ell ) \subset V_\ell  \subset \mathbb R^N \setminus \{ O \}$
for each $\ell $.
Take a homeomorphism
$G_\ell  : \mathbb R^N \cong \mathbb R^N$
such that 
$G_\ell  (B_\ell ) = T(Z_\ell )$ and
$G_\ell (x) = x$ for each $x\in \mathbb R^N\setminus V_\ell $ \cite[Thm. 2.2]{Keldysh}.
As can be seen from
the proof of \cite[Thm. 2.2]{Keldysh},
we may additionally assume that the
restriction of $G_\ell $ onto $B_\ell $ is an orientation-preserving
similarity transformation.
For each $\ell $, take a corrective ``twist''
$\varphi _\ell  : \mathbb R^N \cong \mathbb R^N$
such that
$\varphi _\ell  (x) = x$ for each $x\in\mathbb R^N \setminus T(Z_\ell )$,
and $\varphi _\ell (G_\ell  (K_\ell )) =  T( Y\cap Z_\ell  )$.
Finally, define a homeomorphism $G:\mathbb R^N \cong \mathbb R^N$ by
$$
G  (x) = 
\begin{cases}
 T(O) &\text{for } x=O, 
\\
\varphi _\ell  G_\ell ( x) &\text{for } x\in V_\ell ,\ \  \ell \in\mathbb N,
\\
x &\text{ otherwise. }
\end{cases}
$$
By construction, $G(X) = T(Y)$. This finishes the proof.

\begin{conjecture}
A Krushkal Cantor set $K\subset \mathbb R^N$, $N\geqslant 4$,
does not possess a sequence of Cantor sets
$X_1,X_2,\ldots \subset \mathbb R^N \setminus K$
such that $\{ X_i , i\in\mathbb N \}$ converges 
homeomorphically to $K$, and $(\mathbb R^N , X_i)
\cong (\mathbb R^N , K)$ for each $i$.
\end{conjecture}

\end{document}